\documentclass[12pt]{article}
\usepackage{graphicx}
\usepackage{amsmath,amsfonts,amssymb,amsthm}
\usepackage{bbm,mathrsfs,bm,mathtools}
\title{The copolarity of pseudo-cones}
\author{Rolf Schneider}
\date{}
\date{}
\newcommand{\R}{{\mathbb R}}

\newcommand{\cL}{{\mathcal L}}

\begin{document}
\maketitle

\begin{abstract}
{We show that the copolarity of pseudo-cones has analogous properties as the usual polarity of convex bodies.}\\[2mm]
\end{abstract}

\section{Introduction}\label{sec1}

For a set $K\subseteq\R^n$, the polar set is defined by
$$ K^\circ:= \{u\in\R^n:\langle u,x\rangle\le 1\mbox{ for all }x\in K\}$$
where $\langle\cdot\,,\cdot\rangle$ denotes the scalar product of $\R^n$. Restricted to the family of compact convex sets $K$ containing the origin $o$ in the interior, polarity is an involution, that is, it satisfies $K^{\circ\circ}=K$. Polarity of convex bodies has several interesting and useful properties.

Within their investigation of general classes of dualities, Artstein--Avidan, Sadovsky and Wyczesany \cite{ASW23} have pointed out that to the duality above there exists a dual polarity, which has much less been studied. For $K\subseteq\R^n$, the copolar set is defined by
$$ K^*:= \{u\in\R^n:\langle u,x\rangle\le -1\mbox{ for all }x\in K\}$$
(compared to \cite{ASW23}, we prefer a reflection in the origin). It satisfies $K^{**}=K$ if restricted to the so-called pseudo-cones. A pseudo-cone in $\R^n$ is a nonempty closed convex subset $K\subset \R^n$ satisfying $o\notin K$ and $\lambda K\subseteq K$ for $\lambda\ge 1$. 

The name `copolarity' was suggested by Rashkovskii \cite{Ras17}, while treating a special case. Parallel to \cite[Sect. 5]{ASW23}, but independently, a first detailed study of copolarity (though not under this name) was undertaken by Xu, Li, and Leng \cite{XLL22}, who also coined the name `pseudo-cone' (though with a slightly different definition). Among several other results, they proved that copolarity is essentially the only inclusion reversing involution on the set of pseudo-cones, thus obtaining a counterpart to a characterization of the classical polarity of convex bodies by B\"or\"oczky and Schneider \cite{BS08}.

There are, in fact, more properties of the classical polarity that have counterparts for copolarity. A systematic investigation seems to be missing. We mention that in \cite{XLL22} it has been shown that the suitably defined support function $h_K$ and radial function $\varrho_K$ of pseudo-cones $K$ satisfy 
\begin{equation}\label{1.1}
\varrho_K(x) = \frac{-1}{h_{K^*}(x)}
\end{equation}
for $x$ in the interior of the recession cone of $K$. For a pseudo-cone $K$ whose recession cone is pointed and $n$-dimensional, the exposed faces of $K^*$ in relation to the exposed faces of $K$ were studied in \cite{Sch24}. Here it turned out that one must distinguish whether faces are bounded or unbounded and whether or not they are contained in the boundary of the recession cone.

It is the purpose of the following to collect some more properties of copolarity which are in analogy to known properties of ordinary polarity.

We assume throughout that $K$ is a $C$-pseudo-cone, that is, a pseudo-cone in $\R^n$ with recession cone $C$, where $C$ is a pointed closed convex cone with interior points. Then the dual cone $C^\circ$ has the same properties.

\section{A relation between support functions}\label{sec2}

We define the support function $h_K:\R^n\to\overline\R=\R\cup\{\infty\}$ of $K$ by
$$ h_K(u):= \sup\{\langle u,x\rangle: x\in K\},\quad u\in \R^n.$$
Then $h_K$ is positively homogeneous and a closed convex function (see, e.g., \cite[Thm. 2.7]{HW20}). It is non-positive on $C^\circ$ and negative on ${\rm int}\,C^\circ$. If $u\in\R^n\setminus C^\circ$, then there exists $x\in {\rm int}\,C$ with $\langle u,x\rangle >0$. Since $\lambda x\in K$ for all sufficiently large $\lambda$, it follows that $h_K(u)=\infty$.

We note that, for $u\in \R^n$,
$$ u\in K^*\Leftrightarrow \langle u,x\rangle \le -1\,\forall x\in K \Leftrightarrow \sup\{\langle u,x\rangle:x\in K \}\le -1 \Leftrightarrow h_K(u)\le -1.$$

For a convex function $f:\R^n\to\overline \R$, the Legendre transform is defined by
$$ (\cL f)(u) := \sup\{\langle u,x\rangle -f(x):x\in \R^n\}\quad\mbox{for } u\in \R^n.$$

We would like to apply this to $\textstyle{\frac{1}{2}}h_K^2$, but this function is not convex, and $-\textstyle{\frac{1}{2}}h_K^2$ is only convex on $C^\circ$. Therefore, we define
$$ \widetilde h_K(u):= \left\{\begin{array}{ll} -\textstyle{\frac{1}{2}}h_K^2(u) & \mbox{if } u\in C^\circ,\\[2mm]
\infty & \mbox{if }u\in\R^n\setminus C^\circ.\end{array}\right.$$
Then $\widetilde h_K$ is convex. In a similar way, $\widetilde h_{K^*}$ is defined. With these definitions, we can state that
\begin{equation}\label{2.1n}
\cL\widetilde h_K =  \widetilde h_{K^*}.
\end{equation}
This is analogous to a relation for the ordinary polarity (see, e.g., \cite[(1.48)]{Sch14}), with appropriately changed definitions.

For the proof of (\ref{2.1n}) we note that for $x\in \R^n$ we have
\begin{eqnarray*}
h_{K^*}(x) &=& \sup\{\langle x,u\rangle : u\in K^*\} = \sup\{\langle x,u\rangle: h_K(u)\le -1\}\\
&=&  \sup\{\langle x,u\rangle: h_K(u)= -1\} = \sup_{u\in C^\circ\setminus\{o\}}\{\langle x,u\rangle: h_K(u)= -1\},
\end{eqnarray*}
since $h_K(u)>-1$ if $u\notin C^\circ\setminus\{o\}$. Thus
$$ h_{K^*}(x) = \sup_{u\in C^\circ\setminus\{o\}} \frac{\langle x,u\rangle}{- h_K(u)},\quad x\in \R^n.$$
Here we have used that the function $u\mapsto\langle u,x\rangle/(-h_K(u))$ is positively homogeneous of degree zero. The latter fact is also used below. Observing that
$$ \min_{a\in\R} \{\lambda a+{\textstyle\frac{1}{2}}a^2\} = -{\textstyle\frac{1}{2}} \lambda^2 \quad\mbox{for } \lambda\in\R,$$
we get for $x\in C$ (where $h_{K^*}$ is finite)
\begin{eqnarray*}
\widetilde h_{K^*}(x) &=& -{\textstyle\frac{1}{2}}h^2_{K^*}(x)\\
&=& \min_{a\in\R}\left\{h_{K^*}(x) a+{\textstyle\frac{1}{2}}a^2\right\} = \min_{a\in\R}\left\{\sup_{u\in C^\circ\setminus\{o\}} \frac{\langle x,u\rangle}{- h_K(u)}a+{\textstyle\frac{1}{2}}a^2\right\}\\
&=&  \min_{a\in\R}\sup_{u\in C^\circ\setminus\{o\}}\left\{ \frac{\langle x,u\rangle}{- h_K(u)}a+{\textstyle\frac{1}{2}}a^2: - h_K(u)=a\right\}\\
&=& \sup_{u\in C^\circ \setminus\{o\}}\left\{ \frac{\langle x,u\rangle}{- h_K(u)}a+{\textstyle\frac{1}{2}}a^2: - h_K(u)=a\right\}\\
&=& \sup\{ \langle x,u\rangle + {\textstyle\frac{1}{2}} h_K^2(u):u\in C^\circ\}\\
&=& \sup\{ \langle x,u\rangle -\widetilde h_K(u):u\in\R^n\}\\
&=& \left(\cL \widetilde h_K\right)(x).
\end{eqnarray*}
For $x\in\R^n\setminus C$ there exists $u\in {\rm int}\,C^\circ$ with $\langle x,u\rangle >0$, hence $\lambda u\in K^*$ for all sufficiently large $\lambda$. This gives 
$$ \left(\cL \widetilde h_K\right)(x) = \sup\{\langle x,u\rangle-\widetilde h_K(u):u\in\R^n\} = \infty.$$
Thus, (\ref{2.1n}) is proved.

\section{A relation between gauge functions}\label{sec3}

For a convex body $L\subset\R^n$, containing the origin in the interior, one defines the radial function by $\varrho_L(x)= \sup\{\lambda\ge 0: \lambda x\in L\}$ and the gauge function by $\|x\|_L= \inf\{\lambda\ge 0: x\in\lambda L\}$, thus
$$ \|x\|_L = \frac{1}{\varrho_L(x)} = h_{L^\circ}(x).$$
In calculus of variations and in Finsler geometry (see, e.g., Rund \cite{Run58}), where often $L$ figures as the `indicatrix' and $L^\circ$ as the `figuratrix`, one prefers the notation $\|\cdot\|_L=F$ and $\|\cdot\|_{L^\circ}=H$, and relations between the functions $F$ and $H$ are of interest. In this section, we study analogous functions for $C$-pseudo-cones and their relations.

We assume in this section that the $C$-pseudo-cone $K$ is smooth and strictly convex inside ${\rm int}\,C$. This means that each point $x\in\partial_i K:=\partial K\cap {\rm int}\,C$ is contained in a unique supporting hyperplane of $K$ and that there is a supporting hyperplane $H$ of $K$ with $H\cap K=\{x\}$. We recall from \cite[Def. 4]{Sch24} that a pair $(x,x^*)\in\R^n\times\R^n$ is called a {\em crucial pair} of $K$ if $x\in \partial K$ and $x^*$ is an outer normal vector of $K$ at $x$ such that $\langle x,x^*\rangle =-1$. Every point $x\in \partial_i K$ belongs to some crucial pair of $K$. If $(x,x^*)$ is a crucial pair of $K$ with $x\in{\rm int}\,C$, then $(x^*,x)$ is a crucial pair of $K^*$, by Lemma 6 of \cite{Sch24}. And if $x\in{\rm int}\,C$, then $x^*\in{\rm int}\,C^\circ$, since otherwise $K$ could not be strictly convex. Since we have assumed that $K$ is smooth and strictly convex inside ${\rm int}\,C$, in a crucial pair $(x,x^*)$ of $K$ with $x\in{\rm int}\,C$ each of the vectors determines the other one uniquely. Therefore, the same holds for $K^*$, so that also $K^*$ is smooth and strictly convex inside $C^\circ$. 

In the present case (and writing $\partial_i K^*:= \partial K^*\cap {\rm int}\,C^\circ$), we can define a map
$$ f_K:\partial_i K\to\partial_i K^* \quad\mbox{by}\quad f_K(x):= x^* \quad\mbox{if $(x,x^*)$ is a crucial pair of $K$}.$$
We have $f_{K^*} \circ f_K= {\rm id}$ on $\partial_i K$ and $f_K^{-1}(x^*)=x$. We shall relate this function to the function $F$ defined below.

We define the {\em gauge functions} of $K$ and its copolar pseudo-cone $K^*$ by
$$ F(x):= \max\{\lambda >0:x\in\lambda K\},\quad x\in{\rm int}\,C,$$
$$ H(x^*) := \max\{\lambda>0: x^* \in\lambda K^*\},\quad x^*\in{\rm int}\,C^\circ.$$
We have
$$ K\cap{\rm int}\,C=\{x\in{\rm int}\,C: F(x)\ge 1\}.$$
Defining the radial function $\varrho_K$ of $K$ on ${\rm int}\,C$ by
$$ \varrho_K(x):= \min\{\lambda>0:\lambda x\in K\}, \quad x\in{\rm int}\,C,$$
we have by (\ref{1.1})
$$ F= \frac{1}{\varrho_K} = -h_{K^*}|_{{\rm int}\,C}.$$ 
Since $K^*$ is strictly convex inside ${\rm int}\,C^\circ$, the support function $h_{K^*}$ is differentiable on ${\rm int}\,C$, hence $F$ is differentiable. 

Let $X$ be a $C^1$ parametrization of $\partial_i K$, with local parameters $u_\alpha$. Let $x^1,\dots,x^n$ be the standard coordinates of $X$. Since $F(X)\equiv 1$, we get
$$ 0= \frac{\partial \frac{1}{2}F^2}{\partial u_\alpha}(X) = \sum_{i=1}^n F(X) \frac{\partial F}{\partial x^i}(X) \frac{\partial x^i}{\partial u_\alpha} = \langle \nabla F(X), X_\alpha\rangle,$$
where $\nabla$ denotes the gradient and $X_\alpha$ is the vector with coordinates $\frac{\partial x^i}{\partial u_\alpha}$, $i=1,\dots,n$. Thus, the gradient $\nabla F(X)$ is orthogonal to the tangent vectors of $X$ and hence is proportional to the normal vector. Since $F^2$ is homogeneous of degree $2$, Euler's homogeneity relation gives
\begin{equation}\label{3.0} 
\left\langle\left(\nabla{\textstyle\frac{1}{2}}F^2\right)(x),x\right\rangle = F^2(x), \quad x\in{\rm int}\,C,
\end{equation}
thus
$$ \left\langle -\left(\nabla{\textstyle\frac{1}{2}}F^2\right)(X),X\right\rangle = -1.$$
Hence, for $x\in\partial_i K$,
$$ \left( x, -\left(\nabla{\textstyle\frac{1}{2}}F^2\right)(x)\right)$$
is a crucial pair. It follows that 
$$ f_K(x) = -\nabla{\textstyle\frac{1}{2}}F^2(x) \quad\mbox{for } x\in\partial_i K.$$
Thus, the gradient mapping of $-\frac{1}{2}F^2$ maps $\partial_i K$ to $\partial_i K^*$. Since this mapping is positively homogeneous of degree one, it maps ${\rm int}\,C$ onto ${\rm int}\,C^\circ$.

We state that
\begin{equation}\label{3.2}
 F(x) = H(x^*)\quad\mbox{if $(x,x^*)$ is a crucial pair of $K$}.
\end{equation}

For the proof, we use Remark 1.6.18 of \cite{Sch24} with $f=-\frac{1}{2}F^2$ (which is equal to $-\frac{1}{2}h^2_{K^*}|_{{\rm int}\,C}$ and hence convex), $D={\rm int}\,C$ and $U={\rm int}\,C^\circ$. Let $(x,x^*)$ be a crucial pair of $K$, then $(\nabla f)(x)=x^*$, $(\nabla f)^{-1}(x^*)=x$, and Euler's relation (\ref{3.0}) reads $\langle x,(\nabla f)(x)\rangle = 2f(x)$. We obtain
\begin{equation}\label{3.1} 
(\cL f)(x^*) = \langle (\nabla f)^{-1}(x^*),x^*\rangle - f((\nabla f)^{-1}(x^*)) = \langle x,(\nabla f)(x)\rangle -f(x) = f(x).
\end{equation}
By (\ref{2.1n}) (with $K$ and $K^*$ interchanged) we have $\cL \widetilde h_{K^*} = \widetilde h_K$. Here $\widetilde h_{K^*}|_{{\rm int}\,C} = -\frac{1}{2}h_{K^*}^2|_{{\rm int}\,C} =-\frac{1}{2}F^2$, and similarly $\widetilde h_K|_{{\rm int}\,C^\circ}=-\frac{1}{2}H^2$, hence
$$ \cL(-\textstyle{\frac{1}{2}}F^2)=-\textstyle{\frac{1}{2}}H^2 \quad \mbox{on }C^\circ.$$
Together with (\ref{3.1}) this gives $-\frac{1}{2}H^2(x^*) = \cL(-\frac{1}{2}F^2)(x^*) = f(x)= -\frac{1}{2}F^2(x)$ and hence the assertion (\ref{3.2}).

More can be said if $F$ is twice continuously differentiable. Then we can define
$$G_{ij}(x) := {\textstyle\frac{1}{2}} \frac{\partial^2 F^2}{\partial x^i\partial x^j} (x) \quad\mbox{for }x\in{\rm int}\,C.$$
Since $(\partial/\partial x^i)F^2$ is homogeneous of degree one, the homogeneity relations yield
$$ \sum_{j=1}^n {\textstyle\frac{1}{2}}\frac{\partial^2F^2}{\partial x^i\partial x^j}(x)x^i = {\textstyle\frac{1}{2}} \frac{\partial F^2}{\partial x^i} (x).$$ 
Thus
$$ \sum_{j=1}^n G_{ij}x^i = \mbox{$i$th coordinate of $\nabla\left({\textstyle\frac{1}{2}}F^2\right)(x)$}.$$
In other words, if we denote by $G(x)$ the $n\times n$ matrix with entries $G_{ij}(x)$, then
$$ f_K(x) = -G(x)x\quad \mbox{for }x\in\partial K.$$

\section{A relation between Gauss curvatures}\label{sec4}

In this section we assume that the $C$-pseudo-cone $K$ and its copolar pseudo-cone $K^*$ have sufficiently smooth boundaries with positive Gauss curvatures. Precise assumptions will be given below. For $x\in\partial K$ we denote by $\nu_K(x)$ to outer unit normal vector of $K$ at $x$. If the Gauss curvature of $\partial K$ at the boundary point $x$ is denoted by $\kappa(K,x)$, we shall show that for crucial pairs $(x,x^*)$,
\begin{equation}\label{4.1}
\frac{\langle x,\nu_K(x)\rangle}{\kappa(K,x)^{\frac{1}{n+1}}}\cdot \frac{\langle x^*,\nu_{K^*}(x^*)\rangle}{\kappa(K^*,x^*)^{\frac{1}{n+1}}}=1.
\end{equation}

The term
\begin{equation}\label{4.2} 
\varrho_{\rm aff}(K,x) = \frac{\langle x,\nu_K(x)\rangle}{\kappa(K,x)^{\frac{1}{n+1}}}
\end{equation}
appearing in (\ref{4.1}) expresses in Euclidean terms what in affine differential geometry is known as the equiaffine distance of $\partial K$ from $o$, or as the equiaffine support function of $K$, taken at $x$. A hypersurface of class $C^2_+$ on which the equiaffine support function is constant is known to be an affine hypersphere. Thus we can deduce from (\ref{4.1}) the following observation: If $\partial K$ is an affine hypersphere, then also $\partial K^*$ is an affine hypersphere.
The corresponding result for ordinary polarity of convex bodies can be found, e.g., in  Salkowski \cite[p. 155]{Sal34} (for $n=3$) or Oliker and Simon \cite[p. 107]{OS92}. About (hyperbolic) affine hyperspheres in cones and the solution of corresponding conjectures of Calabi, one can read in Li, Simon and Zhao \cite[§ 2.7]{LSZ93}.  

For convex bodies with a boundary of class $C^3_+$ and the ordinary polarity, relation (\ref{4.1}) can be found for $n=3$ in Salkowski \cite[§ 39, (31), (32)]{Sal34} and for general $n$ in Kaltenbach \cite[Sect. I.3]{Kal90}. We mention that for general convex bodies and suitable sets of crucial pairs (which may be empty, as for example in the case of polytopes), a proof has been given by Hug \cite{Hug96}. We can carry over Kaltenbach's proof to pseudo-cones and copolarity, but since we use a weaker assumption, the proof is slightly more elaborate. Here we use the following assumption.

\noindent{\bf Assumption.} The radial function of the $C$-pseudo-cone $K$ is twice continuously differentiable.

For the proof of (\ref{4.1}), we use local parametrizations of $\partial K$ of the form $X= \varrho_K(U)U$, where $U$ is a regular local parametrization of $\Omega_C$ of class $C^2$, where
$$ \Omega_C:= C\cap{\rm int}\,C$$
is an open subsets of the unit sphere ${\mathbb S}^{n-1}$. Let $N= \nu_K(\varrho_K(U)U)$, so that $N$ is a $C^1$ parametrization of the outer unit normal vector of the parametrized hypersurface $X$. Denoting the coefficients of the first and second fundamental forms of $X$ by $g_{\alpha\beta}= \langle X_\alpha,X_\beta\rangle$, $b_{\alpha\beta} = -\langle N_\alpha, X_\beta\rangle$ and their matrices by $(g_{\alpha\beta})_{\alpha,\beta=1}^{n-1}=g$, $(b_{\alpha\beta})_{\alpha,\beta=1}^{n-1}=b$, we have
$$ \kappa(K,X) = \frac{\det b}{\det g}.$$
From
$$ \frac{\det(N,N_1,\dots,N_{n-1})}{\det(N,X_1,\dots,X_{n-1})}=  \frac{\det(N,N_1,\dots,N_{n-1})\det(N,X_1,\dots,X_{n-1})}{\det(N,X_1,\dots,X_{n-1})\det(N,X_1,\dots,X_{n-1})}= (-1)^{n-1} \frac{\det b}{\det g}$$
we get
\begin{equation}\label{4.4} 
\kappa(K,X) = \left|\frac{\det(N,N_1,\dots,N_{n-1})}{\det(N,X_1,\dots,X_{n-1})}\right|.
\end{equation}

The support function of the copolar pseudo-cone $K^*$ is given by
$$ h_{K^*}(u) = \frac{-1}{\varrho_K(u)} \quad\mbox{for }u\in{\rm int}\, C$$
and hence is of class $C^2$. If $u$ is the exterior unit normal vector of $K^*$ at $x^*$, then
$$ x^*(u)=\nabla h_{K^*}(u),$$
by \cite[Cor. 1.7.23]{Sch14}. Therefore,
$$ X^*= \nabla h_{K^*}(U)$$
gives a local $C^1$ parametrization of $\partial K^*$, with the property that $X^*$ and $X$ have the same local parameters if $(X,X^*)$ is a crucial pair of $K$.

We now use some consequences of the fact that $K^*$ has a support function of class $C^2$ (here we differ from \cite{Kal90}). We apply \cite[pp. 115--117]{Sch14}, noting that the considerations are local, so that the unboundedness of $K^*$ does not make a difference.

First we note that the exterior unit normal vector of $X^*$, denoted by $N^*$, is given by $N^*=U$, hence it is at least of class $C^1$. The correspondence between $K$ and $K^*$ can now be expressed in terms of local parametrizations, namely
\begin{equation}\label{4.5}
N^*=\frac{X}{\|X\|},\quad X^*= \frac{N}{-\langle X,N\rangle},\quad X= \frac{N^*}{-\langle X^*,N^*\rangle},\quad \|X\|=\frac{1}{-\langle X^*,N^*\rangle}.
\end{equation}
We denote differentiation with respect to the $\alpha$th local parameter by an index $\alpha$. By differentiation  we get
\begin{equation}\label{4.6}
N^*_\alpha=\frac{1}{\|X\|}\left( X_\alpha-\frac{\langle X, X_\alpha\rangle}{\|X\|^2}X\right),
\end{equation}
\begin{equation}\label{4.7}
X^*_\alpha = \frac{1}{-\langle X,N\rangle}\left(N_\alpha-\frac{\langle X,N_\alpha\rangle}{-\langle X,N\rangle}N\right).
\end{equation}

From \cite[(2.46)]{Sch14} and the fact that the reciprocal Gauss curvature is the 
 product of the principal radii of curvature (taken, respectively, at a boundary point and the exterior unit normal vector at this point), we obtain that
$$ \kappa(K^*,X^*) = \left|\frac{\det e^*}{\det b^*}\right|,$$
where $e^*$ is the matrix with entries $e^*_{\alpha\beta} = \langle N_\alpha^*,N_\beta^*\rangle$ and $b^*$ has entries $\langle X^*_\alpha,N^*_\beta\rangle$. Similarly as (\ref{4.4}) was obtained, we deduce that
\begin{equation}\label{4.8} 
\kappa(K^*,X^*) = \left|\frac{\det(N^*,N^*_1,\dots,N^*_{n-1})}{\det(N^*,X^*_1,\dots,X^*_{n-1})}\right|.
\end{equation}
This looks like (\ref{4.4}), but we could not just deduce it ``by polarity'', since $\partial K$ is not assumed to be of class $C^3$.

The remaining part follows \cite{Kal90}, with the necessary changes. By (\ref{4.6})
we obtain
\begin{eqnarray*}
&&\det(N^*,N^*_1,\cdots,N^*_{n-1})\\
&&= \|X\|^{-n} \det \left(X, X_1-\frac{\langle X,X_1\rangle}{\|X\|^2}X,\dots,X_{n-1}-\frac{\langle X,X_{n-1}\rangle}{\|X\|^2}X\right)\\
&&= \|X\|^{-n} \det(X,X_1,\dots,X_{n-1}).
\end{eqnarray*}
From this and from 
$$ \frac{\det(X,X_1,\dots,X_{n-1})}{\det(N,X_1,\dots,X_{n-1})} =  \frac{\det(X,X_1,\dots,X_{n-1})\det(N,X_1,\dots,X_{n-1})}{\det(N,X_1,\dots,X_{n-1})\det(N,X_1,\dots,X_{n-1})}= \frac{\langle X,N\rangle}{1}$$
we get
\begin{equation}\label{4.9}
\det(N^*,N^*_1,\cdots,N^*_{n-1})= \|X\|^{-n} \langle X,N\rangle \det(N,X_1,\dots,X_{n-1}).
\end{equation}

In a similar way, we obtain
\begin{eqnarray*}
\det(N^*,X_1^*,\dots,X_{n-1}^*) &=& \frac{1}{\langle X^*,N^*\rangle}\det(X^*,X_1^*,\dots,X_{n-1}^*)\\
&=& \frac{- 1}{\langle X^*,N^*\rangle\langle X,N\rangle}\det(N,X_1^*,\dots,X_{n-1}^*),
\end{eqnarray*}
where (\ref{4.5}) was used. From this and (\ref{4.7}) we get
\begin{equation}\label{4.10} 
\det(N^*,X_1^*,\dots,X_{n-1}^*) =\frac{(-1)^n}{\langle X^*,N^*\rangle\langle X,N\rangle^n}\det(N,N_1,\dots,N_{n-1}).
\end{equation}
Observing that $1/\|X\|=-\langle X^*,N^*\rangle$ by (\ref{4.5}), from (\ref{4.4}), (\ref{4.8}), (\ref{4.9}) and (\ref{4.10}) we obtain
$$ \kappa(K,X)\kappa(K^*,X^*) = \langle X,N\rangle^{n+1}\langle X^*,N^*\rangle^{n+1},$$
which is (\ref{4.1}).

\section{Centro-affine fundamental forms}\label{sec5}

Under stronger smoothness assumptions on $\partial K$ and $\partial K^*$, some information about the fundamental forms of centro-affine differential geometry can be carried over from ordinary polarity to copolarity. For ordinary polarity, these results can be found in Salkowski \cite[§§ 37, 38]{Sal34} (for $n=3$), Laugwitz \cite[III.6]{Lau65}, Oliker and Simon \cite{OS92}. For easier comparison, we use the classical tensor notation employed in \cite{Lau65, Sal34}.

We write crucial pairs of $K$ in the form $(x,x^*)$ (considering $x^*$ as a vector of $\R^n$, since we have tacitly used the scalar product to identify the dual space of $\R^n$ with $\R^n$). We can assume that $\partial K$ is parametrized by parameters $u_\alpha$, $\alpha=1,\dots,n-1$. If the local parametrization of $\partial K$ is given by $X$, then by demanding that $(X,X^*)$ is a crucial pair we have automatically a local parametrization $X^*$ of $\partial K^*$. In other words, we have $X^*= f_K\circ X$, where $f_K$ is the map defined in Section \ref{sec3}. We assume now that both, $X$ and $X^*$, are of class $C^3$ and have positive Gauss curvatures.

The centro-affine differential geometry of a parametrized hypersurface is obtained as its relative differential geometry, if for the relative normal vector one chooses the negative position vector. For the basics of relative differential geometry, we refer to \cite[Sect. 2]{Sch67}, \cite[§ 1]{Sim68}.

We denote differentiation of $X$ or $X^*$ with respect to the local parameter $u_\alpha$ by a lower index $\alpha$. By the definition of the crucial pair $(X,X^*)$ we have
$$ \langle X^*,X_\alpha\rangle =0,\quad \langle X^*,X\rangle = -1.$$
Since $-X$ is the relative normal vector, this shows that $X^*$ is the covector of the tangent plane, as defined in \cite{Sch67}, \cite{Sim68}. Therefore, 
$$ G_{\alpha\beta} := \langle X^*, X_{\alpha\beta}\rangle$$
defines the centro-affine metric tensor of $\partial K$. Denoting by $_{\|\alpha}$ covariant differentiation with respect to the centro-affine Riemannian metric, we also have
$$ G_{\alpha\beta} = \langle X^*, X_{\alpha\|\beta}\rangle = -\langle X^*_\alpha,X_\beta\rangle = - \langle X^*_\beta,X_\alpha\rangle.$$
If $\overline G_{\alpha\beta}$ denotes the centro-affine metric tensor of $\partial K^*$, then we immediately see that
\begin{equation}\label{5.1}
G_{\alpha\beta} = \overline G_{\alpha\beta}.
\end{equation}

The cubic fundamental form of $\partial K$ is defined by 
$$ A_{\alpha\beta\gamma} := \langle X^*, X_{\alpha\|\beta\|\gamma}\rangle.$$
The Ricci relation $G_{\alpha\beta\|\gamma}=0$ gives
$$ \langle X^*,X_{\alpha\|\beta}\rangle_{\|\gamma} =0.$$
hence
$$ \langle X^*_\gamma,X_{\alpha\|\beta}\rangle + \langle X^*, X_{\alpha\|\beta\|\gamma}\rangle =0.$$
Since also $\langle X^*_\gamma,X_\alpha\rangle = -G_{\gamma\alpha}$ and thus $\langle X^*_\gamma,X_\alpha\rangle _{\|\beta} =0$, we get
$$ \langle X^*_{\gamma\|\beta},X_\alpha\rangle + \langle X^*_\gamma,X_{\alpha\|\beta}\rangle =0.$$
The same relation holds with $X^*$ and $X$ interchanged, thus
$$ \langle X^*_{\alpha\|\beta},X_\gamma\rangle + \langle X^*_{\alpha\|\beta\|\gamma}\rangle =0,\quad \langle X^*_\alpha, X_{\gamma\|\beta}\rangle + \langle X^*_{\alpha\|\beta},X_\gamma\rangle =0.$$
Denoting the cubic fundamental form of $\partial K^*$ by $\overline A_{\alpha\beta\gamma}$, we obtain
$$ A_{\alpha\beta\gamma} = -\langle X_\gamma,X_{\alpha\|\beta}\rangle = \langle X^*_{\gamma\|\beta}, X_{\alpha}\rangle = - \overline A_{\gamma\beta\alpha}$$
and thus
\begin{equation}\label{5.2}
A_{\alpha\beta\gamma} = - \overline A_{\alpha\beta\gamma},
\end{equation}
since the cubic fundamental form is symmetric. 

Relations (\ref{5.1}) and (\ref{5.2}) show that the fundamental forms of centro-affine geometry behave under copolarity in the same way as under the ordinary polarity.

\end{document}